\documentclass[preprint,11pt]{elsarticle}

\usepackage{graphicx}
\usepackage{amssymb}
\usepackage{eurosym}
\usepackage{amsmath}
\usepackage{siunitx}
\usepackage{tabularx}
\usepackage{adjustbox}
\usepackage{nomencl}
\usepackage{booktabs}
\usepackage{rotating}
\usepackage{diagbox}
\usepackage{lineno}
\usepackage{enumitem}
\usepackage{etoolbox}
\usepackage{hyperref}
\usepackage[margin=2.5cm]{geometry}

\usepackage{framed} 
\usepackage{multicol} 
\usepackage{multirow} 
\usepackage{nomencl}
\makenomenclature
\usepackage{etoolbox} 

\patchcmd{\pprintMaketitle}
{\hrule\vskip12pt}
{\hrule\vskip12pt\ifvoid\extrainfobox\else\unvbox\extrainfobox\par\vskip12pt\fi}
{}{}

\newsavebox\extrainfobox

\begin{document}

\begin{frontmatter}
	\title{An energy system model for mixed bilateral and pool markets}
	\author[tu]{Ni Wang\corref{cor}}
	\ead{n.wang@tudelft.nl}
	\cortext[cor]{Corresponding author.}
    \author[tu]{Remco A. Verzijlbergh}
    \author[tu]{Petra W. Heijnen}
    \author[tu1]{Paulien M. Herder}
    \address[tu]{Faculty of Technology, Policy and Management, Delft University of Technology, Jaffalaan 5, 2628BX Delft, The Netherlands}
    \address[tu1]{Faculty of Applied Sciences, Delft University of Technology,  Lorentzweg 1, 2628CJ Delft, The Netherlands}

\begin{abstract}
Investments into renewable energy are increasing rapidly around the world. Energy system models are able to provide insights into optimal investment capacities and thus are widely used to aid the long-term investment decision-making under an electricity market environment. Existing energy system models, however, fail to consider bilateral electricity markets while in reality, these constitute a major part of all energy trades. In this paper, we propose an improved energy system model that endogenously considers mixed bilateral and pool markets. In this model, we also introduce three externality cost items that account for the social cost of technologies, carbon taxes/renewable energy subsidies, and the bilateral product differentiation in the bilateral market, respectively. We start with an equilibrium problem formulation for different market players and next, an equivalent optimization problem is presented.
Then, a case study of the pan-European market to reach 95\% emission reduction in 2050 is conducted to demonstrate the model. Different scenarios are constructed to showcase two different usages of product differentiation in the bilateral market, i.e., willingness to pay and exogenous costs. Our main conclusion is that the inclusion of mixed bilateral and pool markets into our enriched energy system model significantly changed the optimal investment capacities, compared to benchmark results from the existing, conventional energy system model. This shows that the inclusion of the bilateral market is of key importance in future investment considerations. Our model is the first of its kind to include this important and realistic bilateral market in energy system models.
\end{abstract}


\end{frontmatter}


\section{Introduction}

\subsection{Background and motivation}
\label{sec:background}
To mitigate climate change and reduce carbon emissions, renewable energy sources (RES) play a vital role in modern power systems. Countries around the world set ambitious RES targets for the next decades and the planning of RES generation and transmission is prominent on the agenda. 

Energy system models refer to optimization models that aim to find the optimal capacity expansion of generation technologies, and transmission networks \cite{Fattahi2020}. The objective is usually to minimize the total system cost while satisfying a number of constraints such as energy balance, generation and transmission limits.
The results of the models are possible scenarios to achieve certain carbon/RES targets that the energy system might evolve into \cite{Pfenninger2014}.
Such models are often used by policy-makers because they serve as a benchmark to help them make decisions on potential policy changes in view of the modeling outcomes, i.e., optimal capacity expansion and the associated costs. For high-RES energy systems, numerous models are built in recent years, see e.g., reviews of \cite{Ringkjob2018, Groissbock2019} on various tools and \cite{Deng2020} for different countries.

Despite the wide use and the policy relevance, existing energy system models have a few characteristics that shift them away from reality. 
First, uncertainty is often not inherently modeled. Long-term investment planning is mainly associated with two types of uncertainties. On the one hand, the parameters used in the models such as cost assumptions and energy demands may be subject to changes in the future \cite{Moret2017}. To account for these uncertainties, energy system models are often complemented with other modeling methods or further analyses such as stochastic programming, robust optimization, sensitivity analysis, and Monte Carlo analysis \cite{DeCarolis2017}.  On the other hand, the formulations of the models may be unable to capture the complexity of the energy system where sub-optimal solutions may be superior, the so-called structural uncertainty \cite{DeCarolis2016}. Here, social factors such as public acceptance and ease of implementation may prevail \cite{Neumann2021}.
Second, the models do not have a comprehensive representation of different electricity markets. 
Classified by time horizon, there are two main types of electricity markets: short-term pool markets that happen on a day before to the day of the energy delivery and long-term bilateral markets that happen from days to years before the actual delivery. 
Energy system models are known to represent the investment equilibrium under a pool market \cite{Ferez-Arriaga1997}, while the bilateral market is not yet included. 
Third, the markets are assumed to be perfect. In electricity markets, to maximize their benefits, producers may not offer the true quantity and the true marginal price to the market, leading to non-socially-optimal results and thus making the market imperfect. 

With regard to the aforementioned features of energy planning which are not covered in energy system models, research efforts have been made that either include those features by modifying the energy system models or different approaches are used. 
Concerning uncertainty, many studies conduct uncertainty analysis or modeling to complement the energy system models, where \cite{Yue2018} gave a systematic review on methods to assess the various uncertainties in energy system models. Meanwhile, market imperfections have been extensively studied in the literature as well. \cite{Gonzalez-Romero2020} reviewed the literature on generation and transmission co-planning under a market environment, where incomplete information and strategic behavior are usually modeled.
However, little attention has been given to the different types of markets in the energy system models. While in practice, a hybrid model which considers the two types of markets would be preferable to the models only considering the separate markets \cite{Lotfi2012}.
Therefore, the motivation of this study is to explore how mixed bilateral and pool markets can be modeled in energy system models.
To that end, we continue to introduce how the bilateral market is currently modeled in the literature.



A bilateral contracting framework is generally known as a bilateral market.
Bilateral contracts are good at hedging against the uncertainties associated with the price and the quantity in short-term markets. 
From the modeling perspective, unlike pool markets, bilateral contracts occur spontaneously and thus are harder to model in a generalized way. 
One potential way to model bilateral contracts in the energy system models is to consider them as a fixed exogenous parameter \cite{Bruno2016}. That is to say, if the bilateral contracts between two parties are already known, the trading volume can be excluded from the energy demands of them.  
This approach is passive, as it can only be used in response to the agreed trading volume, while it does not help the buyers/sellers to determine an optimal trading partner and the associated trading volume. The latter is usually done by optimization models.

The bilateral contracts could also be modeled endogenously. Since the liberalization of the electricity sector, various bilateral market models were proposed to calculate the market equilibrium under different assumptions. These studies focus on producers and their behaviors in the wholesale market. One of the pioneering works in this field is \cite{Hobbs2001}, where Cournot models of imperfect competition were used to simulate the bilateral market. This model was later modified to study generation investment while different carbon policies were evaluated in \cite{He2012}.   
Apart from the equilibrium analysis, research efforts have also been made on individual generator's perspectives to model bilateral contracts. In \cite{El-Khattam2004}, an optimization model was proposed for the optimal planning for distributed generations under competitive market auctions and fixed bilateral contract scenarios.
Other market players than generators such as retailers, prosumers, and energy communities have also been studied. For example, \cite{Karandikar2010} presented a methodology to evaluate bilateral contracts of retailers from a risk perspective.
\cite{Tang2017} proposed a game-theoretical model to describe the competitions for bilateral contracts among generation companies and large consumers.  
\cite{Pourakbari-Kasmaei2020} modeled the trilateral interactions among an integrated community energy system, prosumers, and the wholesale electricity market.
Bilateral contracts were also modeled in combination with demand response to find the optimal energy storage sizing in \cite{Nazari2019}.
In terms of the modeling methods, agent-based modeling is sometimes used to model bilateral contracts. In the review of \cite{Foley2010} on electricity systems models, two agent-based modeling platforms that incorporate bilateral contracts, EMCAS and GTMax, are discussed. 
In addition, \cite{Bower2000} evaluated the effects of bilateral markets in England and Wales using an agent-based simulation. 
\cite{Lopes2012} addressed the challenge of using software agents for the negotiation of bilateral contracts by presenting a multi-agent energy market.
\cite{Imran2020} developed utility-based and adaptive agent-tracking strategies for bilateral negotiations.
Furthermore, \cite{Bompard2008} proposed a complex network approach
for assessing bilateral trading patterns under physical
network constraints.
Lastly, in recent years, with the increasing penetration of distributed energy sources, peer-to-peer (P2P) markets emerge as next-generation market designs.  
In these markets, bilateral trading is considered as one of the most promising P2P market mechanisms \cite{Wang2020c} and is thus commonly modeled. Particularly, bilateral trades can be associated with preferences of the trading parties. To represent this feature, terms such as heterogeneous preferences (\cite{Yang2015, Hahnel2020}), product differentiation \cite{Sorin2019} and energy classes \cite{Morstyn2019a} have been used. Among those, product differentiation is a generic mathematical formulation \cite{Baroche2019a} that can be used for various purposes, e.g., \cite{Baroche2019} used it to account for exogenous network tariffs in P2P markets. 

\subsection{Research gap and contributions}
Based on the background information, we found that despite various methods have been used to model the bilateral market, the hybrid market, i.e., the mixed bilateral and pool markets are not yet modeled in energy system models. 
While energy system models are known for their ease-of-use and policy relevance, the incomplete inclusion of the markets seriously limits the practicality of the models. To overcome this drawback and ensure better inclusion of the markets, in this paper, we propose an improved energy system model that endogenously considers mixed bilateral and pool electricity markets. 
The contributions of this paper are summarized as follows:

\begin{itemize}
    \item We present a partial equilibrium model to represent the long-term investment equilibrium in the mixed bilateral and pool markets for high-RES power systems. The equilibrium model consists of several optimization problems of electricity market participants which will be detailed in Section \ref{sec:concept}.
    \item We propose an improved energy system model where the mixed bilateral and pool markets are modeled endogenously, which is equivalent to the partial equilibrium model. 
    In this model, we introduce three externality cost items that account for the social cost of technologies, carbon taxes/RES subsidies, and the bilateral product differentiation, respectively (see details in \ref{sec:node}). Moreover, apart from the mixed bilateral and pool markets, the proposed planning model is generic so that it can be used to model the bilateral market or the pool market separately.
    \item We conduct a case study for the pan-European market using realistic data. The influences of the mixed bilateral and pool markets on the European power system in 2050 are quantified.
\end{itemize}

The paper is structured as follows. Firstly, Section \ref{sec:concept} provides the preliminaries to understand the model by conceptualizing the mixed bilateral and pool markets in this study. Next, the models are presented in Section \ref{sec:model}. In Section \ref{sec:results}, a case study to illustrate the model is introduced, and then the results are discussed. At last, Section \ref{sec:conclusion} concludes where future recommendations are also given.

\section{Conceptualizing the mixed bilateral and pool markets} \label{sec:concept}

\begin{figure}[htbp]
\centerline{\includegraphics[width=1.3\textwidth, angle=90]{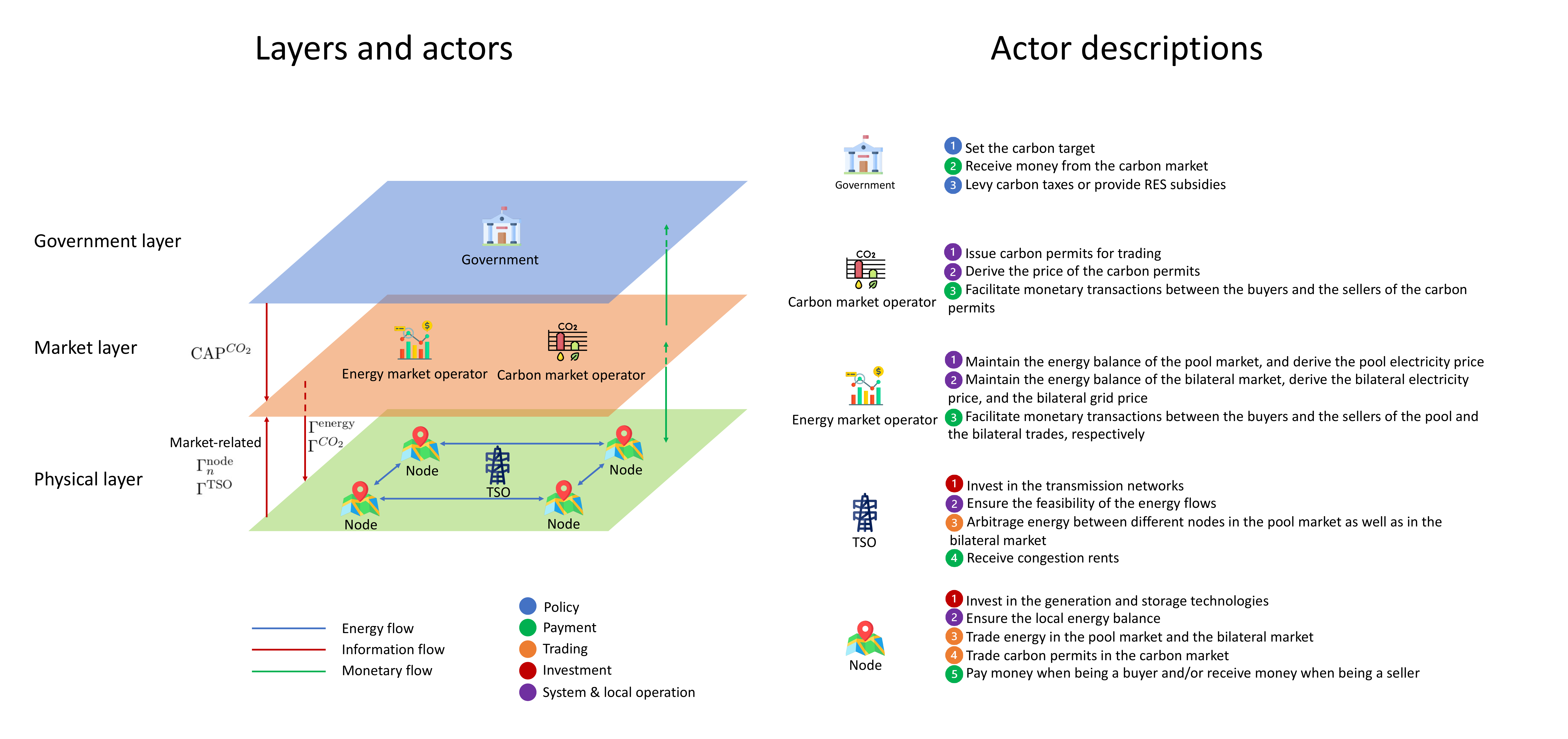}}
\caption{Conceptualization of the mixed bilateral and pool markets: layers, actors and their descriptions.}
\label{fig:concept}
\end{figure}

This section introduces the conceptualization of the mixed bilateral and pool markets in this paper as shown in Figure \ref{fig:concept}. This figure positions the involved actors in several layers and describes the actors by their main responsibilities and activities. The conceptualization will be the basis for the model formulation in Section \ref{sec:model}.

Energy system models are characterized by rich spatio-temporal details of the energy system and accordingly, some simplifications are made to make the model tractable. 
One of the common simplifications is to not explicitly model individual generators and/or consumers, instead, supply and demand that are in proximity are aggregated to geographical nodes. The nodes could be regions in national planning models or countries in continental planning models such as the European planning model. Hence in this study, the nodes are considered as a generic type of actor with an aggregated supply and/or demand, which participates in the market on behalf of the proximity. This way of conceptualizing the problem will be further explained in Section \ref{sec:node}. In addition, the same with the existing energy system models and/or market equilibrium researches such as \cite{StevenA.Gabriel2013}, in this study, the considered actors are limited to a minimal extent while still enough for the functioning of the markets to demonstrate the model. In that regard, some other actors such as retailers or prosumers are left out of scope.    

The lowest layer in Figure \ref{fig:concept} is the physical layer, referring to the generation technologies and the transmission networks with their associated actors.
The capacity expansions of generation and transmission are done by the nodes and the transmission system operator (TSO), respectively.
Regarding energy flows, there are energy exchanges among the nodes, but not between the nodes and the TSO as the TSO neither produces nor consumes energy.

The middle layer is the market layer. There are two market operators, the energy market operator and the carbon market operator. The energy market operator operates the bilateral market and the pool market together. 
In the pool market, locational electricity prices are derived by the market operator, at which the buyers and sellers trade energy. The mechanism is different in the bilateral market.
We utilize the P2P bilateral trading mechanism as proposed by \cite{Baroche2019a}, where nodes negotiate bilaterally about both the amount and the price of the bilateral energy exchange. The market operator helps to facilitate the transactions when a consensus is reached. More information on the bilateral markets will be given in Section \ref{sec:equilibrium}.

Now we turn to the carbon market operator, which will be discussed together with the government level. 
The role of the government is to set carbon goals for the future, which are facilitated through policies such as a cap-and-trade system and carbon taxes/RES subsidies \cite{He2012}.
The cap-and-trade system is a market-based system, where the demands (i.e., the carbon goal) are fixed by the government. A typical example of such markets is the EU Emissions Trading System \cite{EU-ETS}. In addition, the carbon goals could also be implemented indirectly, i.e., through levying carbon taxes or providing RES subsidies. 
To make the model generic, these two policies are both included in this study. However, in practice, one can opt for either or both of the policies depending on the situation.

Finally, there are information and monetary flows between different layers. To introduce these flows, some mathematical symbols are used in Figure \ref{fig:concept}. The mathematical symbols will only be discussed briefly now, since the details will be given in the model formulation in Section \ref{sec:model}.
There is bi-directional information exchange between the actors at the bottom layer and the market operators during the operational phase of the markets.
$\Gamma^{\text{node}}_n$ and $\Gamma^{\text{TSO}}_n$ are the sets of decision variables for node $n$ and the TSO, respectively.
The nodes and the TSO provide the information of energy trades (which are parts of their decision variables, denoted by market-related $\Gamma^{\text{node}}_n$ and $\Gamma^{\text{TSO}}_n$ in Figure \ref{fig:concept}) to the market operators. 
Reversely, the market operators inform them about the prices in the markets, denoted by the set of decision variables of the energy market operator $\Gamma^{\text{energy}}$ and the carbon market operator $\Gamma^{CO_2}$, respectively. 
After the market has been cleared, payments are made and the process is facilitated through the operators as well. On the contrary, the flows between the upper two layers are unidirectional. The government sets a target $\text{CAP}^{CO_2}$ and informs the carbon market operator. The carbon market operator collects the money which goes reversely to the government.

\section{Models}
\label{sec:model}

In this section, the proposed models are presented, which are divided into an equilibrium model that represents the long-term investment equilibrium under mixed bilateral and pool electricity markets, and an improved energy system optimization model that is equivalent to the equilibrium model for the mixed markets.

We start by formulating the respective optimization problems for the nodes, the TSO, the energy market operator, and the carbon market operator. Here, the same as other energy system models, we assume perfect competition and do not address uncertainties inherently in the model. This is because the main purpose of this study is to provide an improved energy system model that incorporates the mixed markets, while other aspects are left out of scope since they have been discussed extensively in the literature as introduced in Section \ref{sec:background}.  
The problem formulations for different actors give a set of optimization problems. Since these problems are all interconnected, i.e., the parameters in one problem may be the decision variables in others, and vice versa, they should be solved together. Essentially, this set of optimization problems forms an equilibrium problem. After obtaining necessary and sufﬁcient conditions for the long-term market equilibrium, the equivalent centralized optimization model will be given. This optimization model is the improved energy system model, since it endogenously models the mixed bilateral and pool markets. The modeling steps, i.e., from the equilibrium model (the collection of individual optimization problems) to the equivalent centralized optimization model, are commonly seen in market equilibrium studies such as \cite{StevenA.Gabriel2013} for the pool market and \cite{Moret2020} for the P2P market.

In this study, lower case symbols are used for variables, upper case symbols are used for parameters. 
Dual variables are expressed using Greek letters and are placed after the colons in the constraints. $n$ is the index for nodes $\mathcal{N}$. $i$ is the index for generation technologies $\mathcal{G}$, and storage technologies $\mathcal{S}$. $l$ is the index for transmission lines in the existing line set $\mathcal{L}$. $t$ represents a time step in the set of total time steps $\mathcal{T}$.  $u$ is a criterion in the criteria set $\mathcal{U}$.

\subsection{Long-term investment equilibrium model}
\label{sec:equilibrium}

\subsubsection{Each node's optimization problem}
\label{sec:node}


Let us first introduce this optimization problem conceptually and then give the mathematical formulations. For all the nodes, their optimization problems will be the same, hence only one node's problem is discussed here.

A node, representing a market participant, wishes to minimize the net cost, or in other words, maximize the net benefit since it can also earn income from sales in the markets. Note that for simplicity, in this study, we assume that a node is a market participant. We are aware that in practice, a node might not always be identical to a market participant, nevertheless, in the same way with the existing energy system models, this model is ready to be extended in case generator/consumer level is of interest, where a node should be further split into more market participants (i.e., generators/consumers). In that respect, we will explain how the model should be modified in the formulation later. 
In this study, the node's optimization problem is cast as a cost-minimization problem. 
The node could take the following actions: investing in generation and storage capacity, producing energy, consuming energy, trading energy in both the bilateral market and the pool market, and trading carbon permits in the carbon market.  
In the bilateral market, a communication graph is pre-defined where the edges connect the pair of the nodes who might trade energy with each other. 
The neighboring nodes on this graph negotiate with each other about the trading volume and the bilateral prices, which is facilitated through the energy market operator.
Note that the communication graph is most of the time not the same as the physical graph where the nodes are connected by electricity networks, since they are different graphs for information exchange and energy exchange, respectively.

$\Gamma^{\text{node}}_n$ is the set of decision variables for node $n$.
It includes the investment capacities $k_{i,n}$ of generation and storage conversion $i$, the investment capacities $k^{\text{storage}}_{i,n}$ of storage $i $, the energy production $p_{i,n,t}$ from technology $i$ at time step $t$, the bilateral trades $p^{\text{bilateral}}_{n,m,t} $ from node $n$ to node $m$ at time step $t$, the pool trades $p^{\text{pool}}_{n,t}$ for node $n$ at time step $t$, state-of-charge $soc_{i,n,t}$ of storage $i$ for node $n$ at time step $t$, storage discharging $p^{\text{out}}_{i,n,t}$ of storage $i$ for node $n$ at time step $t$, storage charging $p^{\text{in}}_{i,n,t}$ of storage $i$ for node $n$ at time step $t$ and the number of carbon permits $e^{CO_2}_n$ for node $n$ to buy from the carbon market in a year.  

The node aims to minimize its total annualized cost related to the investment and the operation of its generation and storage technologies. 
The objective function of the node $n$ is divided into three parts, which are  given in (\ref{equ:obj_non_trading}) - (\ref{equ:obj_externality}).

The first part (\ref{equ:obj_non_trading}) includes the non-trading-related costs, which are Capital Expenditure (CapEx) cost of generation and storage technologies, Fixed Operation \& Maintenance (FOM) costs, Variable Operation \& Maintenance (VOM) costs. These costs are modeled in the same way as in \cite{Wang2020b}. 
Here, $A_i$ is the annuality factor for technology $i$, $C_i$ and $CS_i$ are the CapEx for generation $i$ and storage $i$, respectively, and $B_i$ is the VOM cost for technology $i$.

The second part (\ref{equ:obj_trading}) consists of the trading-related costs including energy trading costs in the pool, energy trading costs and grid costs in the bilateral market, and carbon trading costs.
The node $n$ can trade energy $p^{\text{bilateral}}_{n,m,t} $ at time step $t$ bilaterally with its neighbors $m \in \omega_n$ that are on the communication graph. 
Trading prices $\lambda^{\text{bilateral}}_{n,m,t}$ and grid prices $\lambda^{\text{grid}}_{n,m,t}$ are associated with the energy trades $p^{\text{bilateral}}_{n,m,t}$ in the bilateral market. In the pool market, the node $n$ trades energy $p^{\text{pool}}_{n,t}$ at the price $\lambda^{\text{pool}}_{n,t}$. In the carbon market, the node $n$ buys a certain amount of carbon permits $e^{CO_2}_n$ that are equivalent to their emissions at the carbon price $\lambda^{CO_2}$.


The third part (\ref{equ:obj_externality}) is the externality costs, which consist of three parts. The first two parts refer to the externality costs that are directly related to the generation capacity or the produced energy. One example is the social costs incurred by the social resistance against wind turbines, then $E^{\text{capacity-ex}}_{i,n}$ refers to the unit social cost of wind turbines at node $n$. 
$E^{\text{production-ex}}_{i,n}$ represents either the RES subsidies (resulting in an income) or the carbon taxes (resulting in a cost) for the unit energy produced.
From the modeling perspective, these two terms essentially change the CapEx or the VOM of certain technologies by incorporating these externalities. 
The last externality cost is formulated as a product differentiation term for every bilateral trade referred from \cite{Sorin2019}. 
It is defined as a general cost term related to bilateral trades because the meaning of the costs depends on the interpretation. On the one hand, it may represent exogenous charges related to the bilateral trades, such as transaction costs, tax payments, and network charges. On the other hand, it could be viewed as an improved utility function, representing the willingness to pay for the bilateral trades. These two applications will be further illustrated and discussed in the case study in Section \ref{sec:scenario}. 
Last but not least, the production differentiation coefficient $E^{\text{bilateral-ex}}_{n,m}=\sum_{u \in \mathcal{U}} c^u_n r^u_{n,m}$, where $r^u_{n,m}$ is the trading characteristic under criterion $u \in \mathcal{U}$ and $c^u_n$ is the criterion value for node $n$. This means that the product differentiation term could be a summation of costs for different characteristics, such as distance or a rating that show the popularity of the bilateral trading parties. Further discussions of the product differentiation term can be found in \cite{Sorin2019}. 


\begin{subequations}
\begin{flalign}
	\min_{\Gamma^{\text{node}}_n}
	 & \sum_{i\in (\mathcal{G}+\mathcal{S})} \frac{C_{i} k_{i,n}}{A_{i}} + \sum_{i\in \mathcal{S}} \frac{CS_{i} k^{\text{storage}}_{i,n}}{A_{i}} + \sum_{t \in \mathcal{T}} \sum_{i\in \mathcal{G}}  B_i p_{i,n,t} \label{equ:obj_non_trading} \\ 
	 & + \sum_{t \in \mathcal{T}} \sum_{m \in \omega_n}(\lambda^{\text{bilateral}}_{n,m,t} +  \lambda^{\text{grid}}_{n,m,t})  p^{\text{bilateral}}_{n,m,t}  + \sum_{t \in \mathcal{T}} \lambda^{\text{pool}}_{n,t} p^{\text{pool}}_{n,t} + \lambda^{CO_2} e^{CO_2}_n \label{equ:obj_trading} \\
	 &  + \sum_{i\in \mathcal{G}} E^{\text{capacity-ex}}_{i,n} k_{i,n} + \sum_{i\in \mathcal{G}} \sum_{t \in \mathcal{T}} E^{\text{production-ex}}_{i,n} p_{i,n,t}  + \sum_{t \in \mathcal{T}} \sum_{m \in \omega_n} E^{\text{bilateral-ex}}_{n,m} |p^{\text{bilateral}}_{n,m,t}|  \label{equ:obj_externality}
\end{flalign}  
\begin{flalign}
     \text{subject to: } 
     & \Phi_n (\sum_{i\in \mathcal{G}} p_{i,n,t} - D_{n, t} + \sum_{i\in \mathcal{S}}p^{\text{out}}_{i,n,t} - \sum_{i\in \mathcal{S}}p^{\text{in}}_{i,n,t}) = \sum_{m \in \omega_n} p^{\text{bilateral}}_{n,m,t},  \forall t \in \mathcal{T}: \lambda^{\text{bilateral}}_{n,t}  \label{equ:node_energy_balance1} \\
 	 & (1 - \Phi_n) (\sum_{i\in \mathcal{G}} p_{i,n,t} - D_{n, t} + \sum_{i\in \mathcal{S}}p^{\text{out}}_{i,n,t} - \sum_{i\in \mathcal{S}}p^{\text{in}}_{i,n,t})  = p^{\text{pool}}_{n,t},  \forall t \in \mathcal{T}: \lambda^{\text{pool}}_{n,t} \label{equ:node_energy_balance2} \\  
	 &  0 \le p_{i,n,t} \le E_{i,n,t} (k_{i,n} + K_{i,n})  , \forall i \in \mathcal{G}, \forall t \in \mathcal{T}: \underline \mu_{i,n,t}, \Bar{\mu}_{i,n,t}  \label{equ:node_gen_limit} \\
	 & soc_{i,n,t} = soc_{i,n,t-1} + H^{\text{in}}_{i} p^{\text{in}}_{i,n,t} - \frac{1}{H^{\text{out}}_{i}} p^{\text{out}}_{i,n,t},  \forall i \in \mathcal{S},  \forall t \in \mathcal{T}: \mu^{\text{storage}}_{i,n,t} \label{equ:node_storage_begin} \\
	 & 0 \le soc_{i,n,t} \le k^{\text{storage}}_{i,n} + K^{\text{storage}}_{i,n} , \forall i \in \mathcal{S}, \forall t \in \mathcal{T}: \underline \mu^{\text{storage}}_{i,n,t}, \Bar{\mu}^{\text{storage}}_{i,n,t}  \\
	 & 0 \le p^{\text{out}}_{i,n,t} \le k_{i,n} + K_{i,n}, \forall i \in \mathcal{S},  \forall t \in \mathcal{T}: \underline \mu^{\text{out}}_{i,n,t}, \Bar{\mu}^{\text{out}}_{i,n,t} \\
	 & 0 \le p^{\text{in}}_{i,n,t} \le k_{i,n} + K_{i,n}, \forall i \in \mathcal{S},  \forall t \in \mathcal{T}: \underline \mu^{\text{in}}_{i,n,t}, \Bar{\mu}^{\text{in}}_{i,n,t} \label{equ:node_storage_end} \\
	 & e^{CO_2}_n = W_i \sum_{t \in \mathcal{T}} \sum_{i \in \text{$\mathcal{R}$}} p_{i,n,t}:\lambda^{CO_2}_n   \label{equ:node_CO2}
\end{flalign}
\end{subequations}

(\ref{equ:node_energy_balance1}) and (\ref{equ:node_energy_balance2}) are both the energy balance constraints. 
The net power injection $\sum_{i\in \mathcal{G}} p_{i,n,t} - D_{n, t} + \sum_{i\in \mathcal{S}}(p^{\text{out}}_{i,n,t} - p^{\text{in}}_{i,n,t})$ is divided into two parts: one for the trading in the bilateral market and the other for trading in the pool. On the right-hand side of (\ref{equ:node_energy_balance1}) is the sum of all bilateral trades for node $n$.

$\Phi_n$ is a parameter between 0 - 1 that is determined by the node $n$ itself, indicating the percentage of its net energy that $n$ would like to trade bilaterally, the rest will be traded in the pool.  
This way of modeling the mixed markets strongly aligns with our conceptualization of the markets, that is, one only has to decide ex-ante how much to trade in total in the bilateral market and in the pool market, without determining specifically who to trade with and how much in the bilateral market. Depending on the product differentiation, the model will help the nodes to find the optimal trading partners and the associated trading volumes. In case the trading partners and the associated trading volume are fixed ex-ante, then there are no further decisions to be made and the amount could be deducted from the demands directly.   
Furthermore, this model is generic in that by changing the value of this parameter, the pool market (when $\Phi_n=0$), the bilateral market (when $\Phi_n=1$), or the mixed markets (when $0 < \Phi_n < 1$) can be modeled.

We would also like to continue the discussion on the assumption that the node $n$ is identical to a market participant in this study. As previously introduced, it is a common simplification in energy system models in order to reduce the complexity of the model. In the same way with those models, if desired, one node can further be split into more generators and/or consumers. In that case, both sides of (\ref{equ:node_energy_balance1}) and (\ref{equ:node_energy_balance2}) will be changed such that instead of the net injection and the trades for one node, the sum of the generators and/or consumers will have to be used. Nevertheless, this assumption does not change the gist of the formulation and has no influence on the equivalent centralized optimization model that will be introduced in Section \ref{sec:centralized}. 

(\ref{equ:node_gen_limit}) indicates that the energy production is constrained by the efficiency $E_{i,n,t}$ (capacity factor in case of variable renewable energy) and the capacity of the generation technologies. Here, $K_{i,n}$ is the existing capacity, $k_{i,n}$ is the capacity to be expanded, essentially making the model a capacity expansion model. (\ref{equ:node_storage_begin}) - (\ref{equ:node_storage_end}) are the storage constraints, indicating the change in state-of-charge, and the capacity limits for state-of-charge, charging and discharging, which are modeled the same way as in \cite{Wang2020b}. The last constraint (\ref{equ:node_CO2}) shows that the amount of emissions equals the number of carbon permits.

\subsubsection{TSO's optimization problem}
The role of the TSO is two-fold. First, it ensures the feasibility of the energy flows and accordingly, invests in the transmission network capacity in a cost-optimal manner. 
Second, it harvests congestion rents by trading $z^{\text{bilateral}}_{n,m,t}$ and $z^{\text{pool}}_{n,t}$ in the mixed bilateral and pool markets. 

The decision variables of TSO are represented by the set $\Gamma^{\text{TSO}}$, which includes the investment capacity $k_l$ in line $l$, the bilateral trades $z^{\text{bilateral}}_{n,m,t}$ from $n$ to $m$ at time step $t$, the pool-based trades $z^{\text{pool}}_{n,t}$ for $n$ at time step $t$ and the energy flow $f_{l,t}$ in line $l$ at times step $t$.

The objective function (\ref{equ:tso_obj}) is to minimize the total annualized cost pertaining to its two roles. The first term in (\ref{equ:tso_obj}) is the investment cost for the transmission network, where $\Delta_l$ is the length of line $l$. In addition to the investment cost, the TSO receives the congestion rents from both electricity markets. 

\begin{subequations}
\begin{flalign}
	\min_{\Gamma^{\text{TSO}}} \sum_{l \in \mathcal{L}} \frac{\Delta_l C_l k_l}{A_l} - \sum_{t \in \mathcal{T}} \sum_{n \in \mathcal{N}}( \sum_{m \in \omega_ n}  \lambda^{\text{grid}}_{n,m,t} z^{\text{bilateral}}_{n,m,t} +  \lambda^{\text{pool}}_{n,t}  z^{\text{pool}}_{n, t} ) \label{equ:tso_obj}
\end{flalign}
\begin{flalign}
     \text{subject to: } 
     & f_{l,t} = \sum_{n \in \mathcal{N}} PTDF_{l,n} (\sum_{m \in \omega_n} z^{\text{bilateral}}_{n,m,t} + z^{\text{pool}}_{n, t}),  \forall l \in \mathcal{L},  \forall t \in \mathcal{T}: \lambda^F_{l,t} \label{equ:tso_ptdf} \\
	 & - (k_l + K_l) \le f_{l,t} \le k_l + K_l,  \forall l \in \mathcal{L},  \forall t \in \mathcal{T}: \Bar{\mu}_{l, t}, \underline \mu_{l, t}  \label{equ:tso_limit}
\end{flalign}
\end{subequations}

The energy flow is modeled using direct current power flow equations. In (\ref{equ:tso_ptdf}), the flow $f_{l,t}$ is calculated based on the Power Transfer Distribution Factors (PTDF) matrix and the total net injection at every node $n \in \mathcal{N}$. 
(\ref{equ:tso_limit}) indicates the thermal limits of the energy flows, where $K_l$ is the existing transmission capacity.

\subsubsection{Energy market operator's optimization problem}
The energy market operator clears the mixed markets at each time step $t$ by minimizing the energy imbalances, and thus determines the corresponding prices. 

The set of decision variables $\Gamma^{\text{energy}}$ includes the bilateral trading price $\lambda^{\text{bilateral}}_{n,m,t}$ from $n$ to $m$ at time step $t$, the grid price $\lambda^{\text{grid}}_{n,m,t}$ from $n$ to $m$ at time step $t$, and the pool trading price $\lambda^{\text{pool}}_{n,t}$ for $n$ at time step $t$.

It makes sure that the bilateral trades should be equal in quantity, the trading energy from the node $p^{\text{bilateral}}_{n,m,t} $ is equal to the bilateral arbitraging energy from the TSO $z^{\text{bilateral}}_{n,m,t}$ at each time step $t$, and the pool-based energy trades equal the energy arbitraged by the TSO $z^{\text{pool}}_{n,t}$ at each time step $t$.

\begin{flalign}
	\min_{\Gamma^{\text{energy}}} & \sum_{t \in \mathcal{T}} \sum_{n \in \mathcal{N}} \sum_{m \in \omega_n} \lambda^{\text{bilateral}}_{n,m,t} (p^{\text{bilateral}}_{n,m,t} +  p^{\text{bilateral}}_{m,n,t}) + \sum_{t \in \mathcal{T}} \sum_{n \in \mathcal{N}} \sum_{m \in \omega_n} \lambda^{\text{grid}}_{n,m,t} (p^{\text{bilateral}}_{n,m,t} -z^{\text{bilateral}}_{n,m,t})  \notag \\
	& \sum_{t \in \mathcal{T}} \sum_{n \in \mathcal{N}} \lambda^{\text{pool}}_{n,t} ( p^{\text{pool}}_{n,t}- z^{\text{pool}}_{n, t})
\end{flalign}

\subsubsection{Carbon market operator's optimization problem}
The government can determine the maximum amount of emissions that are allowed to be emitted, which will be regarded as a cap $\text{CAP}^{CO_2}$ in the carbon market. All nodes need to buy the carbon permits from the carbon market which are equivalent to their emissions. The set of decision variables $\Gamma^{CO_2}$, which includes only the carbon price $\lambda^{CO_2}$, will be determined by the carbon market operator. The optimization problem of this operator is formulated as the following. 

\begin{flalign}
	\min_{\Gamma^{CO_2}} \lambda^{CO_2} (\sum_{n \in \mathcal{N}} e^{CO_2}_n - \text{CAP}^{CO_2})
\end{flalign}

\subsubsection{Karush–Kuhn–Tucker conditions}
After laying down the optimization problems of all the actors, it is observed that the decision variables of one problem are the parameters of other problems, i.e., everyone's decisions influence others' decisions. This means that all the optimization problems must be solved together and essentially form an equilibrium problem. The necessary and sufﬁcient optimality conditions, also known as the Karush–Kuhn–Tucker (KKT) conditions, for the long-term investment equilibrium under the mixed bilateral and pool markets are derived. Those conditions are given in \ref{sec:appendix}.

\subsection{Equivalent centralized optimization problem: improved energy system model}
\label{sec:centralized}

In the four optimization problems, the decision variables of one problem only exist in the objective function of others and not in the constraints. This observation indicates that only one solution will exist, which results in a Nash equilibrium.
We follow the approach in market equilibrium studies \cite{StevenA.Gabriel2013, Moret2020} where equivalent optimization problems are derived. We are then able to find the equivalent centralized optimization problem, i.e., the improved energy system model for the mixed bilateral and pool markets, which is formulated as follows as problem (5).

\begin{subequations}
\begin{flalign}
	\min_{\Gamma}
	 & \sum_{i\in (\mathcal{G}+\mathcal{S})} \frac{C_{i} k_{i,n}}{A_{i}} + \sum_{i\in \mathcal{S}} \frac{CS_{i} k^{\text{storage}}_{i,n}}{A_{i}} + \sum_{t \in \mathcal{T}} \sum_{i\in \mathcal{G}}  B_i p_{i,n,t} + + \sum_{l \in \mathcal{L}} \frac{\Delta_l C_l k_l}{A_l} \label{equ:central_obj_first} \\ 
	 &  + \sum_{i\in \mathcal{G}} E^{\text{capacity-ex}}_{i,n} k_{i,n} + \sum_{i\in \mathcal{G}} \sum_{t \in \mathcal{T}} E^{\text{production-ex}}_{i,n} p_{i,n,t}  + \sum_{t \in \mathcal{T}} \sum_{m \in \omega_n} E^{\text{bilateral-ex}}_{n,m} |p^{\text{bilateral}}_{n,m,t}| 
	 \label{equ:central_obj}
\end{flalign}
\begin{flalign}
    \text{subject to: }  
     & (\ref{equ:node_energy_balance1}) - (\ref{equ:node_CO2}), \forall n \in \mathcal{N} \label{equ:central_constraint1} \\
    & (\ref{equ:tso_ptdf}) - (\ref{equ:tso_limit}) \label{equ:central_constraint2}\\
	 & p^{\text{bilateral}}_{n,m,t} = - p^{\text{bilateral}}_{m,n,t}, \forall n \in \mathcal{N}, \forall m \in \omega_n, \forall t\in \mathcal{T}  \text{: } \lambda^{\text{bilateral}}_{n,m,t} \label{equ:central_constraint3} \\
	 & p^{\text{bilateral}}_{n,m,t} = z^{\text{bilateral}}_{n,m,t},  \forall n \in \mathcal{N}, \forall m \in \omega_n,\forall t\in \mathcal{T}   \text{: } \lambda^{\text{grid}}_{n,m,t}  \label{equ:central_constraint4}\\ 
 	 & p^{\text{pool}}_{n,t} = z^{\text{pool}}_{n,t}, \forall n \in \mathcal{N}, \forall t\in \mathcal{T} \text{: } \lambda^{\text{pool}}_{n,t} \label{equ:central_constraint5} \\
	 & \sum_{n \in \mathcal{N}} e^{CO_2}_n = \text{CAP}^{CO_2} \text{: } \lambda^{CO_2} \label{equ:central_co2}
\end{flalign}
\end{subequations}

The decision variables belong to the set $\Gamma$, which includes all the decision variables of the nodes and the TSO. 
The objective function (\ref{equ:central_obj_first}) and (\ref{equ:central_obj}) is the summation of the objective functions of all the actors. 

The constraints are also a gathering of all the constraints of the actors' problems. 
(\ref{equ:central_constraint1}) includes the constraints related to nodal energy balances, generation limits, and storage. (\ref{equ:central_constraint2}) refers to the power flow calculations and the thermal limits of the networks. 
(\ref{equ:central_constraint3}) - (\ref{equ:central_constraint5}) are the KKT conditions of the optimization problem of the energy market operator. 
(\ref{equ:central_constraint3}) is the reciprocity constraint, showing that the bilateral trades should be equal in quantity, where the dual variable $\lambda^{\text{bilateral}}_{n,m,t}$ is the bilateral trading price. (\ref{equ:central_constraint4}) and (\ref{equ:central_constraint5}) are the energy balance constraints between the nodes and the TSO, where the dual variables are the grid price $\lambda^{\text{grid}}_{n,m,t}$ for the bilateral trade and the pool electricity price $\lambda^{\text{pool}}_{n,t}$ for the pool-based trade, respectively. (\ref{equ:central_co2}) gives the cap for all the carbon emissions, with the dual variable $\lambda^{CO_2}$ being the carbon price.

\section{Results and discussions}
\label{sec:results}

To illustrate how the improved energy system model could be used in practice, we explore and show the optimal design of the highly renewable European power system in 2050, where a pan-European electricity market is modeled.
This section presents numerical results from the case study and their implications are discussed.

\subsection{Case study set-up}

\subsubsection{Background}

The EU is committed to combat global climate action under the Paris Agreement and aims to be climate-neutral by 2050 \cite{EU2050}. In this study, we run our model with a 95\% $CO_2$ emission reduction compared to the 1990 level over a year with hourly resolution of weather data and electricity demand for 30 countries. The included generation technologies are offshore wind, onshore wind, solar PV, gas, and hydropower with hydrogen and battery as storage technologies. The 30 countries are connected using high-voltage direct current (HVDC) links. The model assumptions and techno-economic parameters are in line with and could be found in \cite{Schlachtberger2017}. There, they modeled a cost-optimal energy system which essentially indicates a long-term investment equilibrium under a pool market. 

In this study, to further explore the outcomes of the European power system decarbonization under more realistic settings compared to a pool market, the proposed model is utilized to generate results under the mixed bilateral and pool markets. 
In the bilateral market, the countries may choose to trade with others based on their own wills, whereas in the pool market, the energy is traded homogeneously in a pool. 
Here, bilateral contracts are assumed to account for 70 \% of the energy exchanges of each country \cite{EuropeanCommission2021}. A complete communication graph for the 30 countries has been used, i.e., the countries could trade energy bilaterally with each other. The centralized optimization problem (5) is programmed using the Python package Pyomo, and solved using Gurobi solver utilizing the high-performance computational power with a RAM of 160 GB and 16 cores within 2 - 4 hours for different scenarios.   

\subsubsection{Scenario definition}
\label{sec:scenario}

Although all the three externality costs in (\ref{equ:central_obj}) are of significance and interest in practice, we will focus on the last item which is related directly to the bilateral trading and show its effects on the results. From the modeling perspective, the first two items, representing social costs and taxes/subsidies, 
indicate a direct change of cost parameters. The effects resemble global sensitivity study which is a common approach to analyze results in optimization-based researches, and thus are not further investigated in this study. 

Aligning with the two interpretations of the product differentiation term which has been introduced in Section \ref{sec:node}, two corresponding scenarios are set. In addition to a benchmark situation where countries have no preferences in the bilateral market (i.e., $E^{\text{bilateral-ex}}_{n,m}=0$), two other scenarios with different usages of the product differentiation term will be discussed. 

\begin{figure}[htbp]
\centerline{\includegraphics[width=0.75\textwidth]{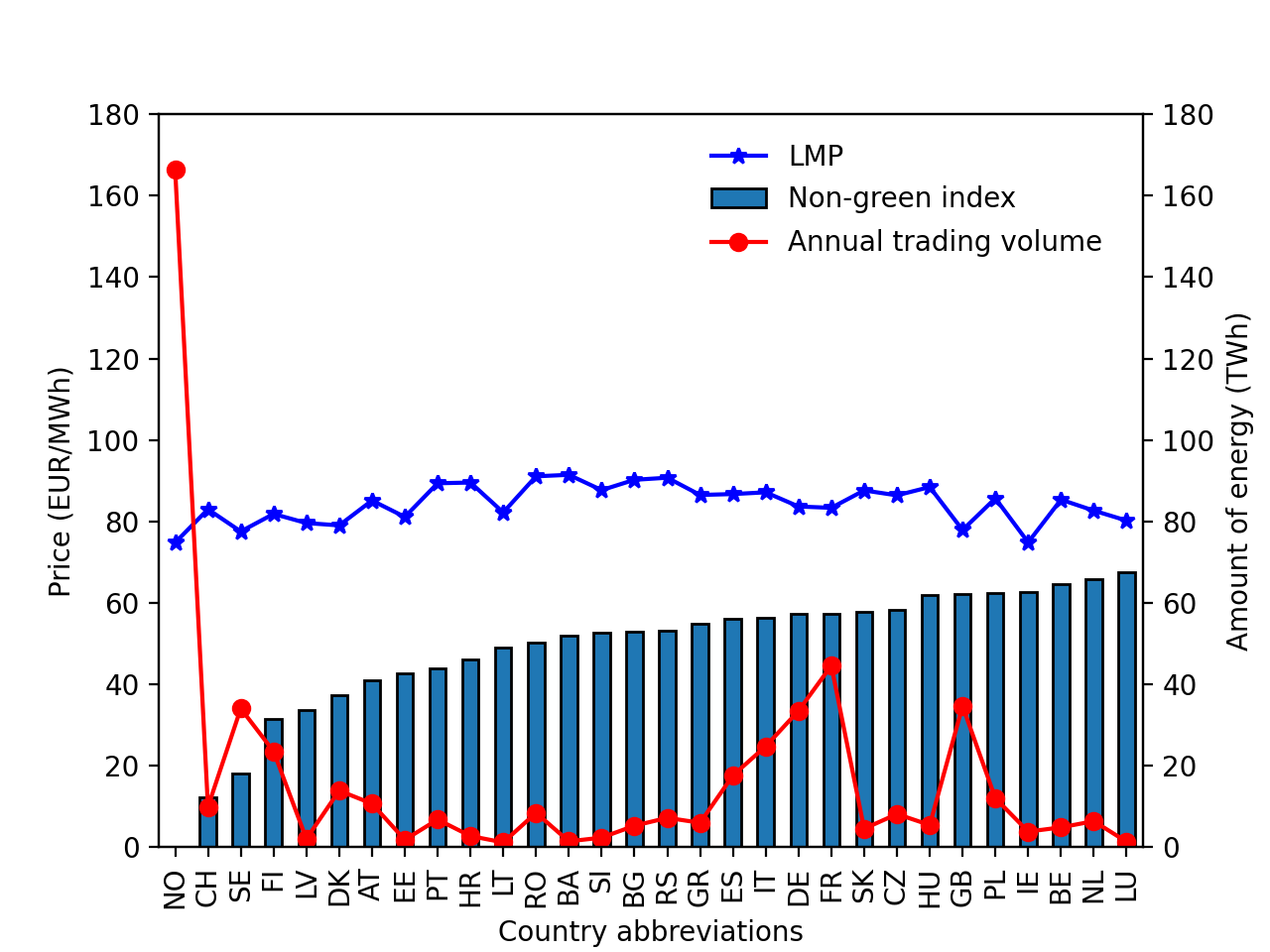}}
\caption{LMP, non-green index, and annual energy trading volume for each country in scenario 1.}
\label{fig:case_III_price}
\end{figure}

In scenario 1, the product differentiation term represents bilateral attitudes (or willingness to pay) between the countries, indicating an improved utility function. This is to show the effects of differentiated preferences on the investment decisions. 
In general, the values of the terms are dependent on the geopolitical relationships. For simplicity, we choose to use only one criterion which is the "non-green index" of each country. 
The rationale is that the product differentiation terms are interpreted as a popularity index among the European countries. 
While it is more realistic if a qualitative study could complement this model to provide data inputs, here we use assumed values for illustrative purposes.
To be more specific, the countries with a higher percentage of RES in their national generation mix are more popular in the bilateral trading for other countries. 
To translate the popularity in costs, we look at the percentage of non-RES generation. For example, Norway is highly dependent on hydropower and only has a non-RES percentage of 25.4 \% in the generation mix. Switzerland ranks second in national RES penetration and has a non-RES share of 37.6 \%. The percentages are translated directly to costs, i.e. 25.4 \EUR{}/MWh and 37.6 \EUR{}/MWh for Norway and Switzerland, respectively. As these costs are derived merely to differentiate the bilateral attitudes, relative values matter more than actual values. We further calibrate all the costs using the lowest one (25.4 \EUR{}/MWh), which ends up with 0 for Norway and 12.2 \EUR{}/MWh for Switzerland. Figure \ref{fig:case_III_price} shows the non-green index of European countries. Luxembourg has the highest unit product differentiation costs (67.6 \EUR{}/MWh) when trading with others.  

\begin{figure*}[htbp]
\centerline{\includegraphics[width=1.2\textwidth]{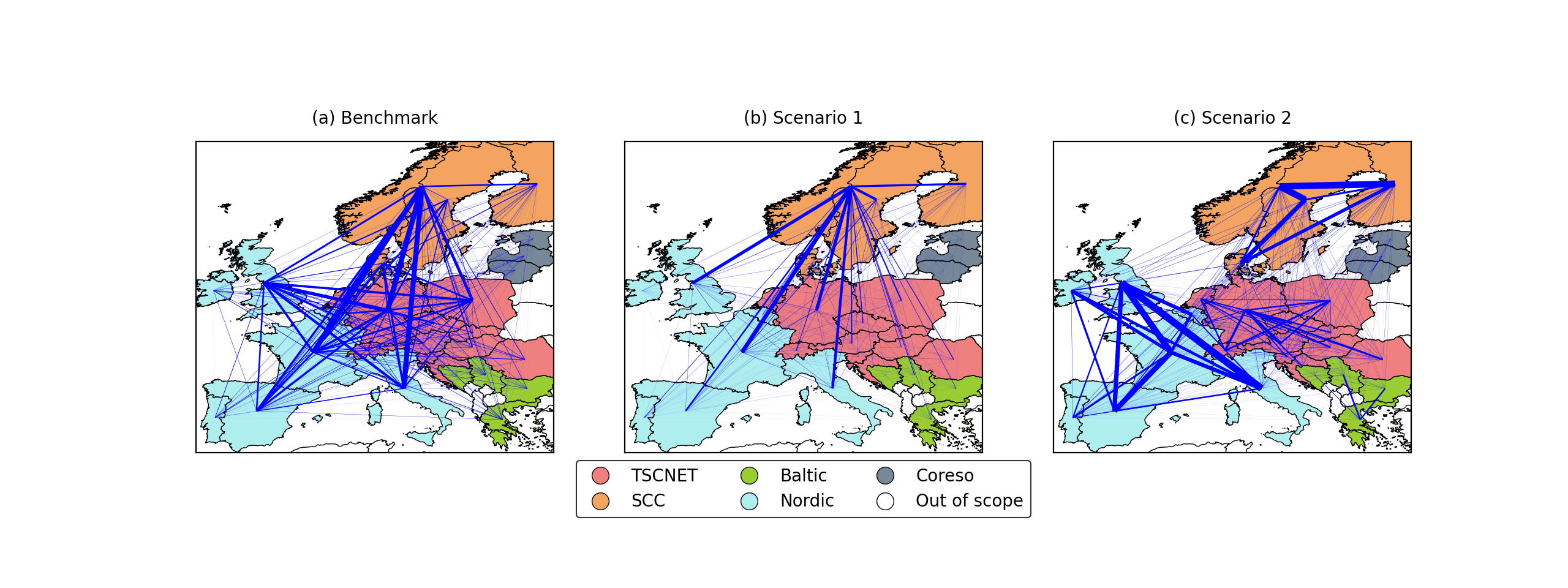}}
\caption{Annual bilateral trading volumes represented by thickness of the lines in benchmark, scenario 1 (with non-green index) and scenario 2 (with inter-regional TC).}
\label{fig:case_III_bilateral}
\end{figure*}

In scenario 2, the terms are used to represent the transaction costs (TCs) between regions. The TCs could be considered as actual costs or deemed as a measure for inter-regional trading barriers, and it is assumed that there are no TCs for intra-regional trades. Different divisions of regions in Europe could be deployed, based on geographical locations or relationships between various multi-national European organizations and agreements. Here, we refer to the TSO regional coordination scheme from ENTSO-E \cite{ENTSO-E_coordination}. In particular, five regional security coordinators are considered, namely TSCNET region, Baltic region, Coreso region, SCC region, and Nordic region, whose geographical locations are shown by different colors in Figure \ref{fig:case_III_bilateral}. 

\subsection{Investment capacities and costs}
In this subsection, we show the results obtained from the different scenarios. First, the total investment capacity of the system will be given. Then, the geographical distribution of the capacities will be shown. Finally, the average costs of countries are discussed. 

\begin{figure}[htbp]
\centerline{\includegraphics[width=0.75\textwidth]{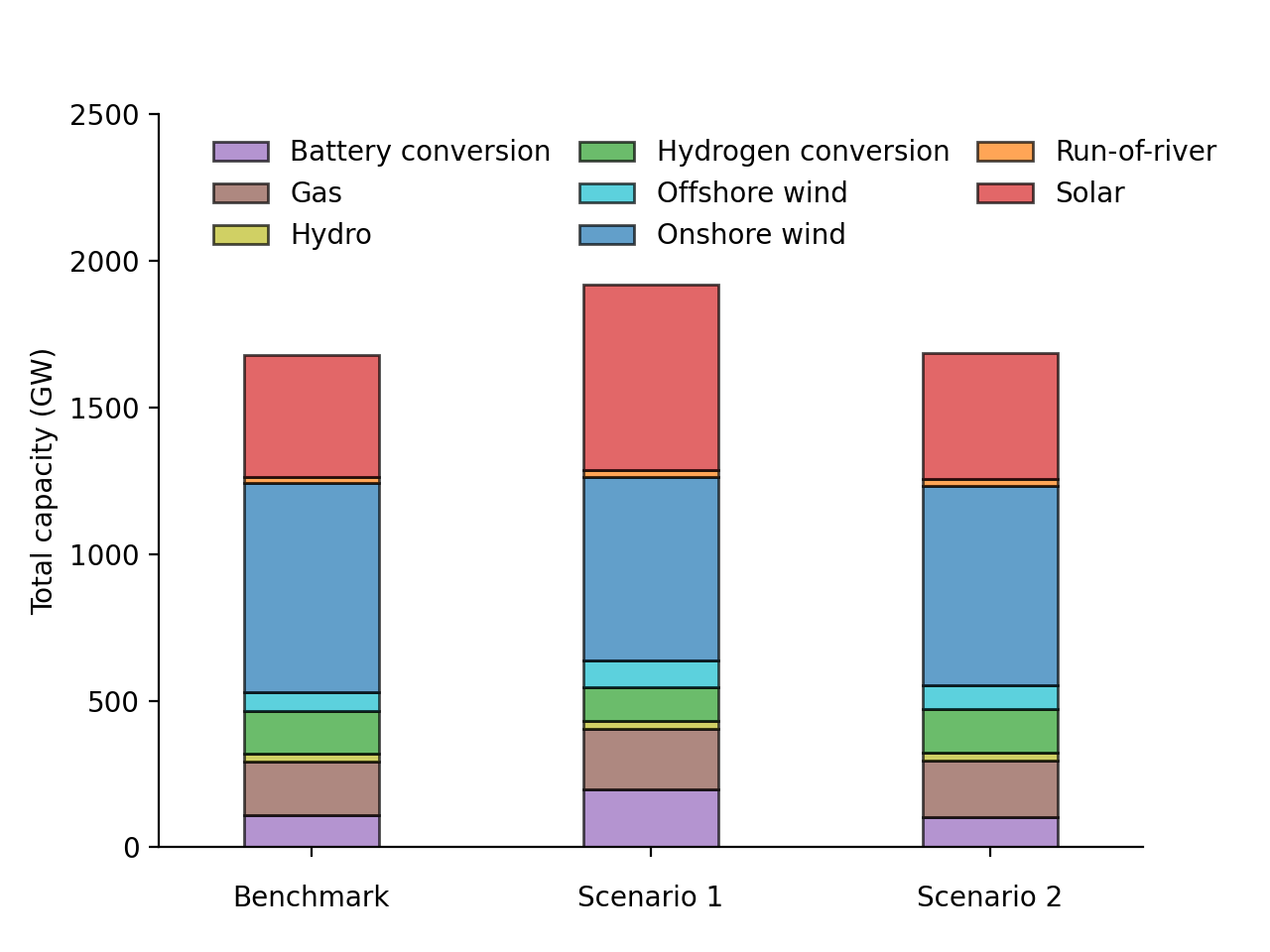}}
\caption{Total investment capacities in the system for benchmark, scenario 1 (with non-green index) and scenario 2 (with inter-regional TC).}
\label{fig:case_III_total_capacity}
\end{figure}

Figure \ref{fig:case_III_total_capacity} illustrates the total investment capacities for the three cases. Several trends could be observed. Firstly, onshore wind energy is the most predominant technology, whose capacity is 42 \%, 33 \%, and 40 \% of the generation mix, respectively in the three cases. Next, solar energy is the second most needed technology and its capacity evens the onshore wind capacity in scenario 1. 
Then, scenario 2 that encourages intra-regional trades shows a similar generation portfolio with that of the benchmark case. Lastly, all other technologies such as offshore wind and storage only contribute a small part to the energy mix.

\begin{figure*}[htbp]
\centerline{\includegraphics[width=1.2\textwidth]{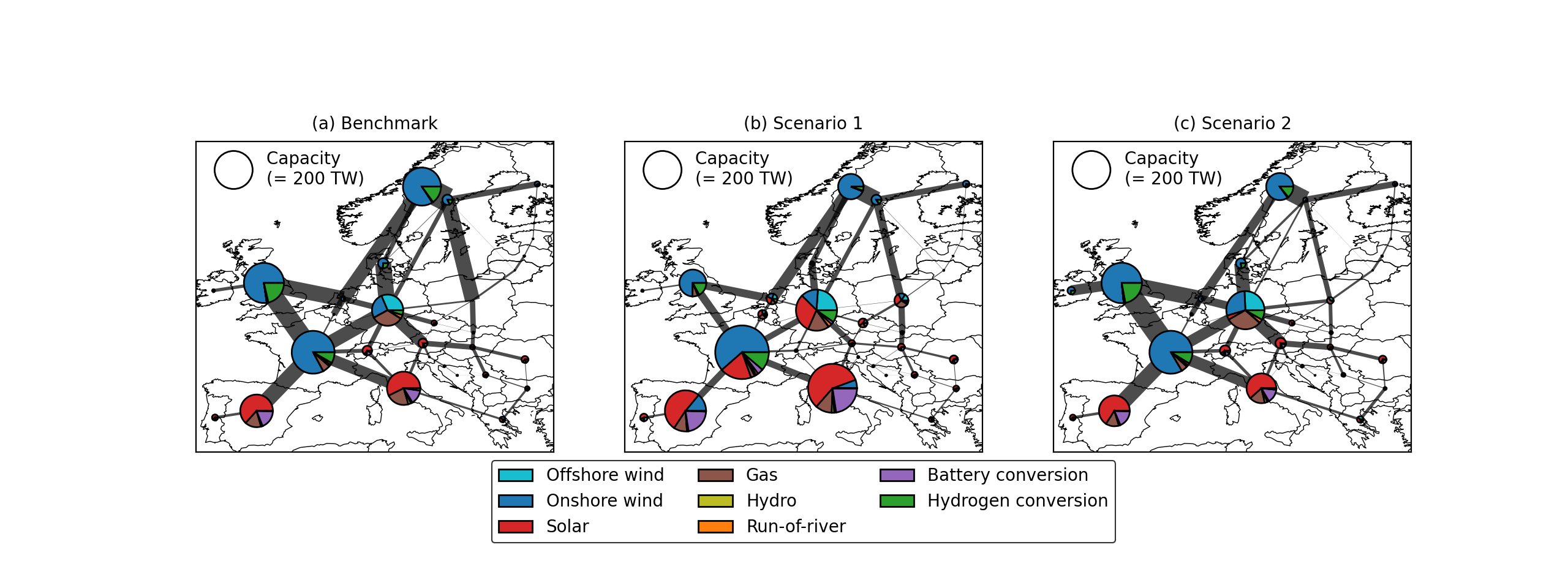}}
\caption{Capacity distribution over Europe for generation, storage and transmission for benchmark, scenario 1 (with non-green index) and scenario 2 (with inter-regional TC). The thickness of the lines represents the transmission capacities.}
\label{fig:case_III_map}
\end{figure*}

We now look at the geographical distribution of the capacities over Europe in Figure \ref{fig:case_III_map}. The meteorological conditions of wind and solar energy are drastically different in various parts of Europe. Therefore, in general, wind capacities are built in the middle and the north, whereas solar capacities are concentrated towards the south. 
Most of the generation capacities are concentrated in the load centers of Europe, such as France and Germany in all the cases. 
In the benchmark scenario where the results are cost-optimal, wind energy capacities are all built in places with good wind resources and transmission networks should be placed around these countries to transmit the energy to other locations.
In particular, Norway builds lots of wind energy capacity and thus becomes a major net export country. 
In scenario 1, the trading positions have changed for many countries due to the introduction of the non-green index, leading to the changed generation and transmission capacities towards southern Europe.
In scenario 2, the geographical capacity distribution is balanced in different regions, indicating more local generations.



\begin{figure*}[htbp]
\centerline{\includegraphics[width=\textwidth]{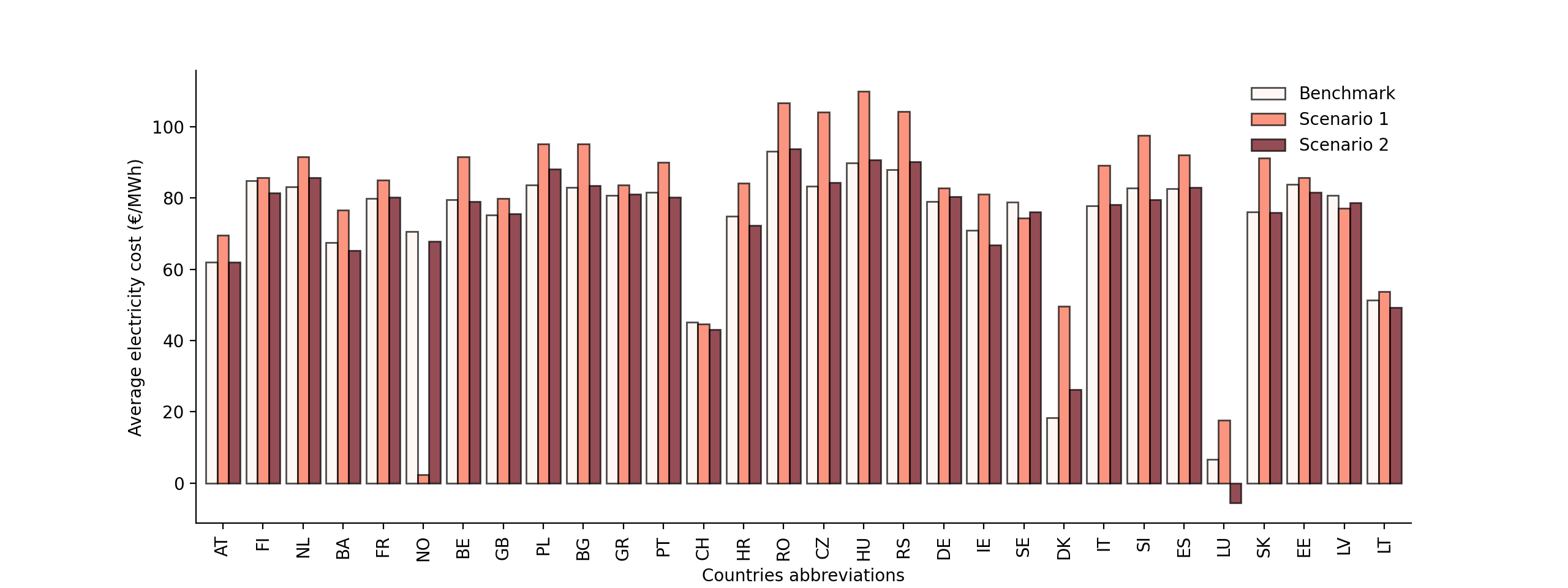}}
\caption{Average costs for the 30 European countries in benchmark, scenario 1 (with non-green index) and scenario 2 (with inter-regional TC).}
\label{fig:case_III_cost}
\end{figure*}

In our framework, the average costs are calculated differently from the traditional cost of electricity. Traditionally, when calculating cost of electricity, only the incurred non-trading-related costs as shown in (\ref{equ:central_obj_first}), i.e., CapEx, FOM and VOM, are included. However, here, we also include the trading costs, which means that the average costs for a country can even be negative provided that lots of revenues are gained from energy trading. Figure \ref{fig:case_III_cost} depicts the average costs per country for the three cases. Overall, the means of the average costs are 73 \EUR{}/MWh, 80 \EUR{}/MWh, 72 \EUR{}/MWh, respectively for the three cases. The costs for most countries stay around these values, with a few exceptions. For example, in scenario 1, Norway is the most popular country to trade with, resulting in the lowest average cost (2 \EUR{}/MWh) therein. Denmark, in contrast, sees a surge in costs in this scenario. This is because it changes from a net exporter to a net importer. It is also found that Luxembourg becomes the only country with a negative average cost in scenario 2. It benefits from the TC-free intra-regional trades and turns its position into a net exporter.  

\subsection{Discussions on the influences of the differentiated non-green index}
After showing the capacities and costs under the three cases, we now take a deeper look at how the differentiated non-green index influences the results in scenario 1. Figure \ref{fig:case_III_price} illustrates the relationships of the non-green index, the locational marginal price (LMP), and the annual trading volume for each country. The bar plots show the non-green index for each bilateral trade. Naturally, one tends to think that the ranking of the non-green index will result in the same ranking for the volumes of the bilateral energy exchanges. While this is true for Norway for being the largest exporter, Switzerland's exports are outnumbered by countries such as France, Germany, Sweden, and Finland. In fact, the volumes of the bilateral energy exchanges are joint effects of several factors, such as the volumes of local energy production in addition to the non-green index. This is the case with countries where the productions are high, such as France and Germany. Therefore, in terms of export percentages, these countries have lower values than Switzerland. In addition to the national energy production, another contributing factor is the price of the bilateral trades. The price not only consists of the product differentiation cost, but also the LMP. In the case of Switzerland, its final price might be higher than other countries, which is another reason to explain its final energy exchange.

\subsection{Discussions on bilateral trades}
Figure \ref{fig:case_III_bilateral} shows the annual bilateral trades for the three cases. In the benchmark case, the bilateral exchanges happen actively across Europe, making use of the different wind and solar resources at different locations. In scenario 1, with various product differentiation costs imposed on the trades, the trading volumes drastically decrease with more locally produced energy. Norway, taking advantage of its lowest non-green index, becomes the most active trader across Europe and is selling energy to all other countries. In scenario 2, trading barriers are created by introducing a TC across regions, which notably change the trading trends. Here, bilateral trades are no longer active across regions, but most of them are kept within regions, showing strong intra-regional preferences.

\section{Conclusions}
\label{sec:conclusion}

Energy system models are known for their policy implications based on the resulted optimal long-term investment decisions, while they have a few key assumptions such as perfect foresight, perfect market, and the pool electricity market environment which limit their applicability in practice.
In this study, we focus on the market environment and present an improved energy system model that endogenously incorporates the mixed bilateral and pool markets. 

The market participants (generalized as nodes in energy system models) invest in generation capacities while taking into account the trades in the market. In this market, they could trade energy with others bilaterally as well as in a pool, which is facilitated through an energy market operator. Moreover, a TSO is responsible for investing in the transmission network and ensures the feasibility of the energy flows in the operational phase. In addition, we utilize a cap-and-trade system for governments to set carbon goals.
The optimization problems of each market participant, i.e., the nodes, the energy market operator, the TSO, and the carbon market operator have been presented. Then, based on the optimality conditions, an equivalent centralized optimization problem has been formulated as the improved energy system model. 
In this model, we introduce three exogenous cost items that account for the social costs of technologies, carbon taxes/RES subsidies, and the bilateral product differentiation, respectively, such that our model includes several exogenous costs generically. While the former two add costs to unit invested capacity or unit produced energy, the product differentiation is imposed on the bilateral trades which could either be viewed as an improved utility function or account for exogenous costs that might be incurred with the bilateral trades, such as transaction costs. 

The model is demonstrated using a case study of the highly renewable European power system in 2050. It has been found that, under the assumed scenarios, the inclusion of the bilateral market significantly changed the results when compared with the conventional cost-optimal energy system model. In terms of the generation mix of the system, the capacity of wind energy drops while that of solar PV increases. The geographical distribution also changes. The cost-optimal results indicate that more generation capacities are placed at locations with favorable weather conditions. However, when bilateral contracts are considered, the results largely depend on other factors, such as the geopolitical relationships between the countries or the multi-national European organizations and agreements. The case study has shown promising results and has demonstrated the urgent need of using this improved energy system model. Our model opens up the possibilities to model future energy systems under more realistic settings considering the mixed markets, while allowing to include several exogenous costs that go beyond the conventional cost-optimal models.       

Nevertheless, the model has a few things that could be improved to warrant future research. On the one hand, the conceptualization of the bilateral contracts in this study is a peer-to-peer model with product differentiation terms. Research efforts could be put into other ways of conceptualizing the bilateral market and/or modeling different market players, and hence the resulting models and insights could be compared with ours. On the other hand, with the trend of sector coupling, it would be interesting to go beyond the power system and investigate how the integration of various energy carriers and even the coupling of energy markets would change the landscape of the European energy system in the future. 

\section*{Acknowledgements}
This research received funding from the Netherlands Organisation for Scientific Research (NWO) [project number: 647.002.007]. Icons in Figure \ref{fig:concept} are made by Freepik from www.flaticon.com.

\bibliographystyle{elsarticle-num-names} 
\bibliography{reference}

\appendix

\renewcommand\nomgroup[1]{%
\vspace{-0pt} 
  \item[\bfseries 
  \ifstrequal{#1}{A}{\emph{Abbreviations}}{%
  \ifstrequal{#1}{S}{\emph{Sets}}{%
  \ifstrequal{#1}{T}{\emph{Superscripts and subscripts}}{%
  \ifstrequal{#1}{V}{\emph{Variables}}{}{%
  \ifstrequal{#1}{P}{\emph{Parameters}}{}}
  }} 
  }%
]}
\newcommand{\nomunit}[1]{%
\renewcommand{\nomentryend}{\hspace*{\fill}#1}}

\vspace{\baselineskip}
\nomenclature[S]{$\mathcal{G}$}{Set of generation technologies}
\nomenclature[S]{$\mathcal{L}$}{Set of lines}
\nomenclature[S]{$\mathcal{N}$}{Set of nodes}
\nomenclature[S]{$\mathcal{R}$}{Set of fossil-fueled generation}
\nomenclature[S]{$\mathcal{S}$}{Set of storage technologies}
\nomenclature[S]{$\mathcal{T}$}{Set of time steps}
\nomenclature[S]{$\mathcal{U}$}{Set of criteria}
\nomenclature[S]{$\omega_n$}{Set of neighbors of node $n$ on the communication graph}
\nomenclature[S]{$\Gamma$}{Set of decision variables}

\nomenclature[P]{$A_{i}$}{Annuity factor of technology $i$}
\nomenclature[P]{$A_l$}{Annuity factor of line $l$}
\nomenclature[P]{$B_{i}$}{Variable Operation \& Maintenance costs (VOM)  of technology $i$ [\EUR{}/MWh]}
\nomenclature[P]{$C_{i}$}{Capital Expenditures (CapEx) and Fixed Operation \& Maintenance costs (FOM) of generation and storage conversion technology $i$ [\EUR{}/MW]}
\nomenclature[P]{$CS_{i}$}{CapEx and FOM of storage technology $i$ [\EUR{}/MWh]}
\nomenclature[P]{$C_l$}{CapEx and FOM of line $l$ [\EUR{}/MW/km]}
\nomenclature[P]{$W_i$}{Emission of fossil-fueled technology $i$ [ton/MWh]}
\nomenclature[P]{$H^{in}_{i}$}{Charging efficiency of storage technology $i$}
\nomenclature[P]{$H^{out}_{i}$}{Discharging efficiency of storage technology $i$}
\nomenclature[P]{$E^{\text{capacity-ex}}_{i,n}$}{Unit social cost for technology $i$ at node $n$ [\EUR{}/MW]}
\nomenclature[P]{$E^{\text{production-ex}}_{i,n}$}{Carbon tax/RES subsidy for technology $i$ at node $n$ [\EUR{}/MWh]}
\nomenclature[P]{$E^{\text{bilateral-ex}}_{n,m}$}{Product differentiation value from node $n$ towards node $m$ [\EUR{}/MWh]}
\nomenclature[P]{$c^u_n$}{Criterion value at node $n$ for criterion $u$}
\nomenclature[P]{$r^u_{n,m}$}{Trading characteristic value from node $n$ towards node $m$ [\EUR{}/MWh]}
\nomenclature[P]{$K_{i,n}$}{Existing capacity of technology $i$ at node $n$ [MW]}
\nomenclature[P]{$K_l$}{Existing capacity of line $l$ [MW]}
\nomenclature[P]{$\Delta_l$}{Length of line $l$}
\nomenclature[P]{$PTDF$}{Power Transfer Distribution Factors matrix}
\nomenclature[P]{$\Phi_n$}{Percentage of trades in the bilateral market at node $n$}
\nomenclature[P]{$\text{CAP}^{CO_2}$}{Carbon market cap}

\nomenclature[V]{$e^{CO_2}_n$}{Carbon permits at node $n$}
\nomenclature[V]{$f_{l,t}$}{Energy flow at line $l$ at time step $t$ [MWh]}
\nomenclature[V]{$k_{i,n}$}{Investment capacity of generation and storage conversion technology $i$ at node $n$ [MW]}
\nomenclature[V]{$k^{\text{storage}}_{i,n}$}{Investment capacity of storage technology $i$ at node $n$ [MWh]}
\nomenclature[V]{$k_l$}{Investment capacity of line $l$ [MW]}
\nomenclature[V]{$p_{i,n,t}$}{Produced energy from technology $i$ from node $n$ at time step $t$ [MWh]}
\nomenclature[V]{$p^{\text{pool}}_{n,t}$}{Pool trades at node $n$ at time step $t$ [MWh]}
\nomenclature[V]{$p^{\text{bilateral}}_{n,m,t}$}{Bilateral trades from node $n$ to node $m$ at time step $t$ [MWh]}
\nomenclature[V]{$p^{\text{in}}_{i,n,t}$}{Charging to storage of storage technology $i$ at node $n$ at time step $t$ [MWh]}
\nomenclature[V]{$p^{\text{out}}_{i,n,t}$}{Discharging from storage of storage technology $i$ at node $n$ at time step $t$ [MWh]}
\nomenclature[V]{$soc_{i,n,t}$}{Energy in storage of storage technology $i$ at node $n$ at time step $t$ [MWh]}
\nomenclature[V]{$z^{\text{bilateral}}_{n,m,t}$}{Arbitrage energy from the transmission system operator (TSO) from node $n$ to node $m$ at time step $t$ in the bilateral market [MWh]}
\nomenclature[V]{$z^{\text{pool}}_{n,t}$}{Arbitrage energy from the TSO at node $n$ at time step $t$ in the pool market [MWh]}
\nomenclature[V]{$\lambda^{\text{bilateral}}_{n,m,t}$}{Bilateral trading price between node $n$ and node $m$ at time step $t$ [\EUR{}/MWh]}
\nomenclature[V]{$\lambda^{\text{grid}}_{n,m,t}$}{Grid price between node $n$ and node $m$ at time step $t$ [\EUR{}/MWh]}
\nomenclature[V]{$\lambda^{\text{pool}}_{n,t}$}{Pool trading price at node $n$ at time step $t$ [\EUR{}/MWh]}

\printnomenclature
\section{Karush–Kuhn–Tucker conditions}
\label{sec:appendix}

This appendix shows the KKT conditions of the equilibrium problem.
Note that the form $0 \le x \perp y \ge 0$ means $x \ge 0$, $y \ge 0$ and $xy = 0$.

{\allowdisplaybreaks

\begin{flalign}
	& \frac{C_{i}}{A_{i}} + E^{\text{capacity-ex}}_{i,n} - \Bar{\mu}_{i,n,t} - \Bar{\mu}^{\text{in}}_{i,n,t} - \Bar{\mu}^{\text{out}}_{i,n,t} = 0,  \forall i\in (\mathcal{G}+\mathcal{S}), \forall n \in \mathcal{N}, \forall t \in \mathcal{T} \\
	& B_i + E^{\text{production-ex}}_{i,n} + \Phi_n \lambda^{\text{bilateral}}_{n,t}  + (1 - \Phi_n) \lambda^{\text{pool}}_{n,t} + \Bar{\mu}_{i,n,t} - \underline \mu_{i,n,t} - W_i \lambda^{CO_2}_n (i \in \mathcal{R})  = 0, \notag \\
	& \forall i\in \mathcal{G}, \forall n \in \mathcal{N}, \forall t \in \mathcal{T} \\
	& \frac{CS_{i}}{A_{i}} - \Bar{\mu}^{\text{storage}}_{i,n,t} = 0, \forall i\in \mathcal{S}, \forall n \in \mathcal{N}, \forall t \in \mathcal{T}   \\
	& \mu^{\text{storage}}_{i,n,t} + \Bar{\mu}^{\text{storage}}_{i,n,t} - \underline \mu^{\text{storage}}_{i,n,t}= 0, \forall i\in \mathcal{S}, \forall n \in \mathcal{N}, \forall t \in \mathcal{T}  \\
	& \frac{\mu^{\text{storage}}_{i,n,t}}{\eta^{out}_{i}} - \Bar{\mu}^{\text{out}}_{i,n,t} - \underline \mu^{\text{out}}_{i,n,t}= 0, \forall i\in \mathcal{S}, \forall n \in \mathcal{N}, \forall t \in \mathcal{T} \\
	& \mu^{\text{storage}}_{i,n,t} \eta^{in}_{i} - \Bar{\mu}^{\text{in}}_{i,n,t} - \underline \mu^{\text{in}}_{i,n,t} = 0, \forall i \in \mathcal{S}, \forall n \in \mathcal{N}, \forall t \in \mathcal{T}  \\
	& E^{\text{bilateral-ex}}_{n,m} \text{sign}(p^{\text{bilateral}}_{n,m,t})  + \lambda^{\text{bilateral}}_{n,t} + \lambda^{\text{grid}}_{n,m,t} + \lambda^{\text{bilateral}}_{n,m,t} = 0,  \forall n \in \mathcal{N}, \forall m \in \omega_n, \forall t \in \mathcal{T} \\
	& \lambda^{CO_2} - \lambda^{CO_2}_n = 0,  \forall n \in \mathcal{N} \\	
    & \frac{\Delta_l C_l}{A_l} - \Bar{\mu}_{l, t} - \underline \mu_{l, t}  = 0,  \forall l \in \mathcal{L} \\
	& -\lambda^{\text{grid}}_{n,m,t} - \lambda^F_{l,t} \sum_l \text{PTDF}_{l,n} = 0,  \forall n \in \mathcal{N}, \forall m \in \omega_n, \forall t \in \mathcal{T} \\
	& -\lambda^{\text{pool}}_{n,t} - \lambda^F_{l,t} \sum_l \text{PTDF}_{l,n} = 0,  \forall n \in \mathcal{N}, \forall t \in \mathcal{T} \\
	& \lambda^F_{l,t} +  \Bar{\mu}_{l, t} - \underline \mu_{l, t} = 0,   \forall l \in \mathcal{L},  \forall t \in \mathcal{T}\\
    & \Phi_n (\sum_{i\in \mathcal{G}} p_{i,n,t} - D_{n, t} + \sum_{i\in \mathcal{S}}p^{\text{out}}_{i,n,t} - \sum_{i\in \mathcal{S}}p^{\text{in}}_{i,n,t}) = \sum_{m \in \omega_n} p^{\text{bilateral}}_{n,m,t},  \forall t \in \mathcal{T},  \forall n \in \mathcal{N} \\
 	& (1 - \Phi_n) (\sum_{i\in \mathcal{G}} p_{i,n,t} - D_{n, t} + \sum_{i\in \mathcal{S}}p^{\text{out}}_{i,n,t} - \sum_{i\in \mathcal{S}}p^{\text{in}}_{i,n,t})  = p^{\text{pool}}_{n,t},  \forall t \in \mathcal{T},  \forall n \in \mathcal{N}\\
 	& soc_{i,n,t} = soc_{i,n,t-1} + H^{\text{in}}_{i} p^{\text{in}}_{i,n,t} - \frac{1}{H^{\text{out}}_{i}} p^{\text{out}}_{i,n,t},  \forall i \in \mathcal{S},  \forall t \in \mathcal{T},  \forall n \in \mathcal{N} \\
	& e^{CO_2}_n = W_i \sum_{t \in \mathcal{T}} \sum_{i \in \text{$\mathcal{R}$}} p_{i,n,t}, \forall n \in \mathcal{N} \\
    & f_{l,t} = \sum_{n \in \mathcal{N}} PTDF_{l,n} (\sum_{m \in \omega_n} z^{\text{bilateral}}_{n,m,t} + z^{\text{pool}}_{n, t}),  \forall l \in \mathcal{L},  \forall t \in \mathcal{T} \\
	& p^{\text{bilateral}}_{n,m,t} = - p^{\text{bilateral}}_{m,n,t}, \forall n \in \mathcal{N}, \forall m \in \omega_n, \forall t\in \mathcal{T} \\
	& p^{\text{bilateral}}_{n,m,t} = z^{\text{bilateral}}_{n,m,t},  \forall n \in \mathcal{N}, \forall m \in \omega_n,\forall t\in \mathcal{T} \\
 	& p^{\text{pool}}_{n,t} = z^{\text{pool}}_{n,t}, \forall n \in \mathcal{N}, \forall t\in \mathcal{T} \\
	& \sum_{n \in \mathcal{N}} e^{CO_2}_n = \text{CAP}^{CO_2}\\
    & 0  \le p_{i,n,t} \perp \underline \mu_{i,n,t} \ge 0, \forall i\in \mathcal{G}, \forall n \in \mathcal{N}, \forall t \in \mathcal{T}\\
    & 0  \le E_{i,n,t} (k_{i,n} + K_{i,n}) - p_{i,n,t}  \perp \Bar{\mu}_{i,n,t} \ge 0, \forall i\in \mathcal{G}, \forall n \in \mathcal{N}, \forall t \in \mathcal{T}\\
	& 0 \le soc_{i,n,t} \perp \underline \mu^{\text{storage}}_{i,n,t} \ge 0, \forall i\in \mathcal{S}, \forall n \in \mathcal{N}, \forall t \in \mathcal{T}\\    
	& 0 \le k^{\text{storage}}_{i,n} + K^{\text{storage}}_{i,n} - soc_{i,n,t} \perp \Bar{\mu}^{\text{storage}}_{i,n,t} \ge 0, \forall i\in \mathcal{S}, \forall n \in \mathcal{N}, \forall t \in \mathcal{T}  \\
	& 0 \le p^{\text{out}}_{i,n,t} \perp \underline \mu^{\text{out}}_{i,n,t} \ge 0, \forall i\in \mathcal{S}, \forall n \in \mathcal{N}, \forall t \in \mathcal{T}\\
	& 0 \le k_{i,n} + K_{i,n} - p^{\text{out}}_{i,n,t} \perp  \Bar{\mu}^{\text{out}}_{i,n,t}  \ge 0, \forall i\in \mathcal{S}, \forall n \in \mathcal{N}, \forall t \in \mathcal{T}  \\
	& 0 \le p^{\text{in}}_{i,n,t} \perp \underline \mu^{\text{in}}_{i,n,t} \ge 0, \forall i\in \mathcal{S}, \forall n \in \mathcal{N}, \forall t \in \mathcal{T} \\
	& 0 \le k_{i,n} + K_{i,n} - p^{\text{in}}_{i,n,t}  \perp  \Bar{\mu}^{\text{in}}_{i,n,t}  \ge 0, \forall i\in \mathcal{S}, \forall n \in \mathcal{N}, \forall t \in \mathcal{T}  \\
    & 0 \le f_{l, t} + k_l + K_l  \perp \underline \mu_{l, t} \ge 0,  \forall l \in \mathcal{L}, \forall t \in \mathcal{T}  \\
    & 0 \le k_l + K_l - f_{l, t} \perp \Bar{\mu}_{l, t} \ge 0,  \forall l \in \mathcal{L}, \forall t \in \mathcal{T}
\end{flalign}

}

\end{document}